\documentclass
[12pt]{article}
\usepackage{amssymb}
\newtheorem{proposition}{\quad\large Proposition}
\newtheorem{lemma}{\quad\large Lemma}
\newtheorem{corollary}{\quad\large Corollary}
\begin{document}
\title{Bernoulli Sieve\\
}
\author{Alexander V. Gnedin\\
{\it Department of Mathematics, Utrecht University}\\
{\it Postbus 80010, 3508 TA Utrecht, The Netherlands}}
\date{}
\maketitle

{\bf Abstract.} Bernoulli sieve is a recursive construction of a random composition (ordered partition)
of integer $n$.
This composition can be induced by sampling from a random discrete distribution which has frequencies
equal to the sizes of component intervals of a 
stick-breaking interval partition of $[0,1]$.
We exploit Markov property  of the composition and its renewal representation to derive asymptotics of the moments
and to prove a central limit theorem for the number of parts.

\vskip0.5cm

\par {\bf 1.}  The Bernoulli sieve  can be seen as a generalisation of
the  `game' found in \cite{Bruss}.
The first round of the game
 starts with 
$n$ players and amounts to tossing a coin with probability  $X_1$ for tails. Each of the players  tosses one time
and the players flipping tails must drop out. If all $n$ get heads the trial is disqualified and must be
repeated completely with all $n$ players, as many times as necessary until some players do quit.
If at least one player remains after the first round, the second round continues with the remaining players,
who must toss another coin with probability $X_2$ for tails.
The game lasts with probabilities $X_3,X_4,\ldots$ for tails until all players are sorted out.
It is assumed  that the probabilities  $X_1,X_2,\ldots$ are independent random variables with a given distribution 
$\omega$ on $]0,1[\,$,
and that given $X_j$ the individual outcomes at round $j$ are conditionally independent. It follows readily
that, as far as only the number of players is concerned,
the outcome of a round depends on the past solely through the number of players which proceed that far.
\par A random composition $C_n$ of integer $n$ arises, with part $j$ being the number of players
dropping out at round $j$. In this paper we shall focus on some properties of $C_n$, in particular we are interested
in the distribution
of the number of parts of the composition, which may be thought of
as the duration of the game.

\par There is a natural way to settle all ${C}_n$'s on the same probability space in a consistent fashion.
Consider a random interval partition of $[0,1]$ by points
$1-(1-X_1)(1-X_2)\ldots (1-X_j)$, $j=1,2,\ldots$ and assign each player a random uniform tag, independent
of the $X_j$'s.
The tags group within  the intervals, and
recording the cluster  sizes, from the left to the right, yields
a composition (intervals containing no tags are ignored). To establish equivalence with 
the coin-tossing construction we only need to note that the chance for a particular player to remain for at least $j$ rounds 
in the game is precisely $(1-X_1)(1-X_2)\ldots (1-X_j)$. 

\par The game in \cite{Bruss} corresponds to 
$\omega$ supported by a single point, in which case
the $X_j$'s  are all equal. The composition is induced then
by sampling 
from a geometric distribution and, of course, had appeared many times in the literature under different guises. 
Karlin  
distinguished this case in the context of a general occupation problem with infinitely many boxes 
and derived distributions for the number of parts, number of singletons, doubletons etc. The feature studied in \cite{Bruss}, \cite{Kirsch} was the probability 
that there is exactly one winner - a player  remaining at the last round (that is to say, the last
part of composition is 1).

\par When $\omega$ is a Beta(1,$\theta$) distribution
the law of $C_n$ is known as the
ordered Ewens sampling formula (ESF). This structure is well understood,
see \cite{ABT} for a recent account and \cite{CSP}, \cite{TSF} for generalisations.

\par Our interest in the construction arose in 
connection with the regenerative compositions 
\cite{GP}.  Within this more general setting the Bernoulli sieve composition 
may be seen as discretisation of a subordinator with  finite   L{\'e}vy measure and zero drift. 
In what follows we shall  treat general measures,
with the only constraint that $\omega$ is not supported by a geometric sequence like $(1-x^{j})$ (in particular,
 sampling from
the geometric distribution 
is ruled out) and such that $\omega$ does not settle too much mass near the endpoints of $[0,1]$.
Our method relies on  renewal theory and the analysis of `divide-and-conquer' recurrences,
the techniques intended to replace the independence-based tools available in  the ESF case, see \cite{ABT}.

\vskip0.5cm

\par {\bf 2.} By exhangeability among the players the
compositions $C_n$  are {\it sampling consistent} for different values of $n$.  That is to say, 
if a part of $C_n$ is selected at random, in a size-biased fashion, and decremented by one unit
then the resulting composition
of $n-1$ (possibly with fewer parts) has the same distribution as $C_{n-1}$.
The sequence $(C_n)$ forms a composition structure in the sense of \cite{RCS},\cite{JMVA} and determine a random
exchangeable composition of a countable set.
\par  There are two further  constructions of $C_n$ featuring renewal and Markov properties.

\vskip0.5cm
\par {\it The  renewal representation} is obtained from 
the stick-breaking construction by 
applying   
transformation $\phi(x)=-\log (1-x)$
which maps $[0,1]$ onto $[0,\infty]$. Consider the range ${\cal R}$  of a renewal process with initial state 0 and
step distribution
$\Omega=\omega_{\phi}$, and let $E_1,\ldots,E_n$
be increasing order statistics from the standard exponential distribution (which correspond to exponentially
distributed tags).
The points of ${\cal R}$ induce a partition of  $[0,\infty]$ into intervals making up the compliment 
${\cal R}^c=[0,\infty]\setminus {\cal R}$,
 and
the points $E_j$ group within the intervals; in these terms composition
$C_n$ becomes a record of all nonzero cluster sizes, from the left to the right.

 \vskip0.5cm
\par {\it Markov chain representation} of $C_n$ stems from the following 
first-part deletion property of $C_n$.
Given the first part of $C_n$ is $m$, the composition of $n-m$ obtained by removing this part has the same distribution
as $C_{n-m}$. This property is obvious in the renewal context: it follows from the regenerative property of ${\cal R}$ (applied
at the leftmost point of ${\cal R}$ to the right of $E_1$) taken together with  the memoryless property of 
the exponential distribution. The deletion property implies that the parts of $C_n$ can be viewed as decrements of a decreasing Markov chain
$Q_n$ which has state-space $\{0,1,\ldots,n\}$, starts at state $n$ and eventually gets absorbed at $0$.
The  one-step transition probability from $n$ to $n-m$ is 
\begin{equation}\label{qw}
q(n,m)=\frac{w(n,m)}{1-w(n,0)},\qquad m=1,\ldots,n.
\end{equation}
where
$$w(n,m)={n\choose m} \int_0^1 x^m (1-x)^{n-m} \omega (\,{\rm d}x),\qquad m=0,1,\ldots,n$$
and the binomial moments of $\omega$. Similar expression can be given in terms of $\Omega$,
with $1-e^{-z}$ in place of $x$.
The quantity $1-w(n,0)=w(n,1)+\ldots +w(n,n)$ will appear throughout as a normalising factor, we  use therefore
the shorter notation $W(n)=1-w(n,0)$. In other terms 
$$
W(n)=\int_0^{\infty} (1-e^{-nz})\,\Omega({\rm d}z)
$$
is the characteristic exponent of the measure $\Omega$ thought  of as a L{\'e}vy measure associated with ${\cal R}$.

\vskip0.5cm

\par  For a given composition $(n_1,\ldots,n_k)$ of  $n$ the probability that $C_n$ assumes this value
is of the product form
\begin{equation}\label{pq}
p(n_1,\ldots,n_k)=q(n_1+\ldots +n_k,n_2+\ldots + n_k)q(n_2+\ldots+n_k,n_3+\ldots +n_k)\cdots q(n_k,n_k),
\end{equation}
because this is  the probability that the chain $Q_n$ has decrements $n_1,\ldots,n_k$
before absorption in 0.

\vskip0.5cm

\par {\bf 3.} We will be interested
in the first instance in the number of parts $K_n$ of composition $C_n$. 
It follows from  $\omega\{1\}=0$   that  $q(n,m)>0$ for all $n\geq m\geq 1$ and $K_n$ goes to infinity
with $n$.

\par Observe that the sizes of intervals comprising the  partition of $[0,1]$  are 
$Y_j=(1-X_1)\ldots (1-X_{j-1})X_j.$ 
Rephrasing
the stick-breaking interpretation, $K_n$ is the number of boxes occupied by at least one of $n$ balls, with probability 
$Y_j$ of hitting the $j$th box.
 Karlin's paper \cite{Karlin} is a basic reference on the model with 
infinitely many boxes and
nonrandom  frequencies, and some information on $K_n$ can be extracted from Karlin's results by conditioning on
$(Y_j)$.

\par Consider two conditions on $\omega$ which limit concentration of mass near $1$ and $0$
\begin{eqnarray}\label{log-b}
\mu:=\int_0^1 |\log(1-x)|\,\omega({\rm d}x)<\infty\,,
\\ \label{llog}
\int_0^1 |\log x|\,\omega({\rm d}x)<\infty\,.
\end{eqnarray}
Reformulated, the condition (\ref{log-b}) says that the first moment of $\Omega$ is finite:
\begin{equation}\label{rho}
\mu=\int_0^{\infty}z\,\Omega({\rm d}z)<\infty\,.
\end{equation}

\par The unconditional law of large numbers from \cite{Karlin} implies

\begin{proposition} 
If conditions {\rm(\ref{log-b})} and {\rm(\ref{llog})} are satisfied then, as $n\to\infty$,
$$
K_n\sim \frac{1}{\mu}\,\log n 
$$
with probability one.
\end{proposition}

{\it Proof.} By the strong law of large numbers we have for $j\to\infty$
$$
-\frac{1}{j}\log Y_j= -\frac{1}{j}\sum_{i=1}^{j-1}\log (1-X_j)-\frac{1}{j}\log X_j\to \frac{1}{\mu},
$$
(condition (\ref{llog}) is necessary and sufficient to have the second term negligible). From this relation we
have (for  Karlin's alpha  on p. 376 of \cite{Karlin})
$$
{\#\{j:Y_j>1/x\}}\sim  \frac{1}{\mu}\,\log x\,,\qquad{\rm as\,\,\,\,} x\to\infty
$$
almost surely. By \cite{Karlin}, Theorem $1'$, 
$$
{\mathbb E}(K_n|(Y_j))\sim \frac{1}{\mu}\,\log n\,,
$$
and by Theorem 8 from that paper the statement holds conditionally on $(Y_j)$, hence also unconditionally.
$\Box$ 

\vskip0.5cm

\par Note that `deconditioning' itself does not allow to conclude about the asymptotics for ${\mathbb E}K_n$ 
(see Proposition 2 to follow).
Results of \cite{Karlin} could be used further to derive the asymptotics of the
conditional variance of $K_n$ and to obtain a conditional central limit theorem. We will not dwell on converting
these results into their unconditional counterparts, rather will take an approach based on the renewal
features of our model.
\vskip0.5cm

\par {\bf 4.} Let $F_n$ be the first part of $C_n$, with distribution ${\mathbb P}(F_n=m)=q(n,m)$.
Markov property of the composition implies that $K_n$ 
satisfies a distributional equation
\begin{equation}\label{K}
K_n\stackrel{d}{=} 1+ K'_{n-F_n}
\end{equation}
where $F_n, K'_1, K'_2,\ldots $  are independent and each $K'_j$ has same distribution as $K_j$. 
Averaging in (\ref{K}) we see 
that $a_n={\mathbb E}\, K_n$    satisfies a linear recursion 
\begin{equation}\label{a}
a_n=1+\sum_{m=1}^{n} q(n,m) a_{n-m}\,
\end{equation}
with boundary value $a_0=0$.
\vskip0.5cm
\par {\large Remark.} Recursions akin to (\ref{a})  
are common in the average-case analysis of algorithms, see references in \cite{RR}.
A recent
dissertation by Bruhn \cite{Bruhn} 
is devoted solely to them. Some results of Bruhn are reproduced in 
R{\"o}sler \cite{Roesler} along with 
distributional analysis of equations more general than  (\ref{K}). 
The class of recursions treated in the cited work relates to the assumption that 
the weights
$q(n,\cdot)$, considered as measures  with support $\{1/n,2/n,\ldots,n/n\}$,
 satisfy an equiboundedness
condition and converge weakly to some measure on $[0,1]$.
\par In our case the convergence of $q(n,\cdot)$ to $\omega$ is clear from  the convergence of moments
(which amounts to Bernstein's trick used to prove the Weierstrass uniform approximation theorem).
Above that, the Bruhn-R{\"o}sler conditions certainly hold when $\omega$ has a smooth density.
However, we were unable to 
check their (very technical) conditions for the general measures $\omega$ and will
rely on the special structure (\ref{qw}). 
A specific feature  of the class of recursions studied here is that
we have a canonical renewal process as a part of the model, while Bruhn and R{\"o}sler needed
to construct an auxiliary renewal process to `mimic' the recursion.

\vskip0.5cm

\par {\large  Remark.} We  formulate the next fact as
 $L^1([0,1],\omega)$-approximability of the logarithm by the (ordinary) Bernstein polynomials,
but in the  sequel we will also make use of the formula (\ref{firstsum}).
There is a variety of closely  related results  in the literature: the 
best known is the aforementioned argument due to Bernstein,
 then there is
a number of $L^1$-results 
on  generalised Bernstein polynomials 
\cite{Lorenz}, and  pointwise  asymptotic expansions found in 
\cite{Flaj-Bern} and \cite{Jacquet}. Still, summation formulas (\ref{firstsum}), (\ref{secondsum}) seem to be new.

\vskip0.5cm
\par The Bernstein polynomial of degree $n$ for $\log (1-x)$ is 
$$
{\rm B}_n(x)=\sum_{m=1}^{n-1}{n\choose m}x^{m}(1-x)^{n-m} \log\left(1-\frac{m}{n}\right)\,.
$$

\vskip0.5cm

\begin{lemma} If $\omega$ satisfies {\rm (\ref{log-b})} then 
$$\lim_{n\to\infty}\int_0^1 |{\rm B}_n(x)-\log(1-x)|\,\omega( {\rm d}x)=0\,.$$
\end{lemma}
{\it Proof.}
There is no simple formula for the 
expectation of the logarithm of  binomial  random variable but
 replacing the logarithms by the harmonic numbers, as $\log (1-m/n)=h_{n-m}-h_n+o((n-m)^{-1})$,
we have the explicit summation formula:
\begin{equation}\label{firstsum}
\sum_{m=0}^{n-1}{n\choose m}x^{m}(1-x)^{n-m}(h_{n-m}-h_n)-x^nh_n=-\frac{x}{1}-\frac{x^2}{2}-\ldots -
\frac{x^n}{n}\,.
\end{equation}
By monotone convergence the series in the RHS approaches
$\log (1-x)$  in the sense of $L^1(\omega,[0,1])$ 
whatever $\omega$.
Getting back to ${\rm B}_n$ easily yields the claim.$\Box$

\vskip0.5cm

\par Proposition 1 strongly suggests the logarithmic asymptotics for $a_n$.
Our proof of this fact will rely on the following simple observation.
Given $n_0>1$ suppose $(a_n)$  satisfies (\ref{a}) for $n\geq n_0$, then
  $(a_n+c)$  also satisfies the recursion for $n\geq n_0$, whatever  constant $c$.

\vskip0.5cm

\begin{proposition} If $\omega$ meets {\rm (\ref{log-b})} then any sequence $(a_n)$ satisfying 
{\rm (\ref{a})} for $n\geq n_0\geq 1$ has asymptotics
$$
a_n\sim \frac{\log n}{\mu}\,.
$$
In particular, this holds for the sequence $a_n={\mathbb E}K_n$ which is the unique solution
which satisfies {\rm (\ref{a})} for $n>0$ and has the boundary 
value $a_0=0.$

\end{proposition}
{\it Proof.} Assume that there exists $\epsilon>0$ such that 
$a_n>(1+\epsilon)\mu^{-1}\log n$
for infinitely many values of $n$. We will lead this to contradiction. 
Selecting $\epsilon$ smaller, for any fixed $c$ we could have inequality 
$a_n>(1+\epsilon)\mu^{-1}\log n +c$
for infinitely many values of $n$. Let $n(c)$ be the minimum such $n$,
then $n(c)\to \infty$ as $c\to\infty$.
Thus for $n<n(c)$ we have 
$
a_n<(1+\epsilon)\mu^{-1}\log n+c
$
which implies 
$$
1+\sum_{m=1}^{n(c)}q(n(c),m)a_{n(c)-m}<1+c+\frac{(1+\epsilon)}{\mu} \,\sum_{m=1}^{n(c)}q(n(c),m)
\log (n(c)-m).
$$
Now from (\ref{a}) and the definition of $n(c)$ we derive
\begin{equation}\label{inter}
(1+\epsilon)\frac{\log n(c)}{\mu} +c<1+c+\frac{(1+\epsilon)}{\mu} \,\sum_{m=1}^{n(c)}q(n(c),m)\log (n(c)-m)
\end{equation}
where $c$ itself cancels but
 $n(c)$ can be taken arbitrarily large by the choice of $c$. 
\par
From Lemma 1 we see that
$$\sum_{m=1}^{n-1} q(n,m)
\log(n-m)=\log n -\mu +o(1)$$
and substituting this formula into (\ref{inter}) and letting $c\to\infty$ yields $0<-\epsilon\,$, which is
 the promised contradiction.
Thus the assumption was wrong and because $\epsilon$ was arbitrary we have
$$\lim\,\sup \frac{a_n}{\mu^{-1}\log n}\leq 1.$$
\par A symmetric argument proves the analogous lower bound, and the claim follows.$\Box$

\vskip0.5cm

\par Turning to the variance of the number of parts $v_n={\rm Var}\, K_n$ 
 we derive from (\ref{K}) a recursion
\begin{equation}\label{v}
v_n=\left(2a_n-1-a_n^2+\sum_{m=1}^n q(n,m) a_{n-m}^2 \right)+\sum_{m=1}^n  q(n,m) v_{n-m}\,,\qquad v_0=0
\end{equation}
which involves $a_n={\mathbb E}K_n$. 
Both (\ref{a}) and (\ref{v}) are instances of the general equation
\begin{equation}\label{br}
b_n=r_n+\sum_{m=1}^{n} q(n,m) b_{n-m}, \qquad b_0=0,
\end{equation}
where $(b_n)$ are unknowns, and $(r_n)$ is given. The proof of Proposition 2 is easily extended to obtain

\vskip0.5cm
\begin{corollary} Assume {\rm (\ref{log-b})}. For any $n_0$ and $r\neq 0$, 
if $(b_n)$ satisfies {\rm (\ref{br})} for $n>n_0$ and if $r_n\to r$ then $b_n\sim r\mu^{-1}\log n$ as $n\to\infty$.
\end{corollary}

\vskip0.5cm

\par With  
a logarithmic asymptotics for $v_n$ in mind, we aim to show the convergence of the bracketed
inhomogeneous term in (\ref{v}). It is easily seen that for this purpose we need 
more than just the principal-term
asymptotics of the expectation, and it is exactly the point where the renewal theory provides indispensable
tools.

\vskip0.5cm

\par {\bf 5.} It is well known that a renewal process starting at $0$ admits a  delayed version
which has the expected number of renewals  within $[0,z]$ (the potential measure) growing linearly with $z$,
see \cite{Feller}. It turns out that  the stationary renewal process induces a `stationary' version of the
Markov chain $Q_n$, which can be used for the asymptotic analysis of (\ref{K}).

\par Let $g(n,m)$ be the probability that $Q_n$ ever visits state $m$ (which means
that at some round of the game there are exactly $m$ players left).
Since $Q_n$ can visit each nonabsorbing state at most once
$g(n,m)$ is also the potential function, i.e. the expected number of visits to $m$.
Interpreting  $r_m$ as a `reward' collected at visit  to state $m$,  we can 
think of $b_n$ satisfying (\ref{br}) as the total
expected reward of $Q_n$. The interpretation implies
\begin{equation}\label{b}
b_n=\sum_{m=1}^{n-1} g(n,m)r_m
\end{equation}
and reduces solving (\ref{a}) to computation of  the potential function. 
Explicit formula is complicated as it involves summation of products 
over a constrained set of compositions of $n$.  Fortunately, 
there is a simple asymptotic formula.

\par Suppose $\Omega$ is not supported by a lattice, and has finite first moment (\ref{rho}).
For $\omega$ this means (\ref{log-b}) and that 
the support is not a 
geometric sequence like $1-x^j$
(in particular, the case of geometric frequencies, when $\omega$ is supported by a single point, is excluded).
Switching to the renewal representation, we introduce a probability distribution 
$$\Omega_0 [0,z]=  \frac{1}{\mu}\int_0^z \Omega[\zeta,\infty]\,{\rm d} \zeta.$$
Let the {\it overshoot} 
$B(z)$ be the distance from $z$ to the leftmost point of ${\cal R}$ to the right of $z$
($B(z)$ is sometimes called the forward process, or forward recurrence time,  or residual lifetime etc.).
The renewal theory, as presented in vol. 2 of the Feller's textbook, says that $\Omega_0$ is the limiting 
distribution of the overshoot as $z\to\infty$.
Observe that $Q_n$ visits $m$ when there is a point 
of $ {\cal R}$ between $E_{n-m-1}$ and $E_{n-m}$ or, equivalently, when the overshoot at $E_{n-m-1}$
does not exceed  $E_{n-m}-E_{n-m-1}$.
The spacing between the two order statistics is independent of $E_{n-m-1}$  and its distribution is Exponential($m$).
By the renewal theorem the distribution of $B(E_{n-m-1})$ converges to $\Omega_0$ as $n\to\infty$ because
$E_{n-m-1}\to\infty$ (in probability), 
thus
$$
g(n,m)={\mathbb P}(B(E_{n-m-1})<E_{n-m}-E_{n-m-1})\to \int_0^{\infty} e^{-mz}\,\Omega_0({\rm d}z)=
\frac{1}{\mu\,m} \int_0^{\infty} (1-e^{-mz}) \,\Omega({\rm d}z),
$$
where the last step follows via integrating by parts.
Changing measure back to $\omega$ we obtain
\vskip0.5cm

\begin{proposition}   If $\omega$  is not supported by a geometric sequence and satisfies {\rm (\ref{log-b})}
then for any $m$ 
$$
\lim_{n\to\infty} g(n,m)=\frac{W(m)}{\mu\, m}\,.
$$
\end{proposition}

\par The proposition suggests to modify chain $Q_n$ 
so that the potential function  becomes exactly 
$$
g_0(m):=\frac{W(m)}{\mu\,m}\,, \qquad
m=1,\ldots,n-1.$$
We shall do this by assuming a special distribution for the first transition
(which can be thought of as a qualifying  round before the game).

\vskip0.5cm 
\par {\large Remark.} Another possibility were to introduce a proper initial distribution on $\{0,1,\ldots,n\}$
so that the formula for potential function were valid  also for $m=n$. But this would correspond 
to composition of a random integer, a model we wish to avoid.
\vskip0.5cm
\par Renewal theory offers construction of a stationary 
version of $ {\cal R}$.
Take $Z_0$ independent of ${\cal R}$  and with distribution $\Omega_0$.
The shifted set ${\cal R}_0=Z_0+{\cal R}$ is the 
range of
the stationary (delayed) renewal process.
For any $z\geq 0$ the overshoot distribution for ${\cal R}_0$ at 
$z$ coincides with $\Omega_0$.

\par
The points of ${\cal R}_0$ induce an interval partition of $[0,\infty]$ thus also a partition of
the sequence of  order statistics $E_1,\ldots,E_n$. Recording the sizes  of blocks
 we obtain a
{\it stationary} composition $C_{0n}$ of $n$.
The parts of ${C}_{0n}$ are considered as  decrements of a new  Markov chain 
${Q}_{0n}$.
Repeating the argument which lead us to Proposition 3 we derive from  invariance of the distribution of $B(z)$
that $g_0$ is the  potential function  of ${Q}_{0n}$.

\par For any reward function the solution of (\ref{br}) satisfies
\begin{equation}\label{ga}
\sum_{m=1}^{n-1} {g}_0(m)r_m = \sum_{m=1}^{n-1} {q}_0(n,m) b_{n-m},
\end{equation}
where ${q}_0(n,\cdot)$ is the distribution of the first part of ${C}_{0n}$.
This formula follows by computing the total expected reward of ${Q}_{0n}$ upon departure
from state $n$. Including state $n$ leads to 
\begin{equation}\label{ga+n}
r_n+\sum_{m=1}^{n-1} {g}_0(m)r_m = \sum_{m=0}^n {w}_0(n,m) a_{n-m},
\end{equation}
where ${w}_0(n,\cdot)$ is the distribution of 
 the number of 
$E_j$'s to the left of $Z_0$.
Explicitly,
$$
{w}_0(n,m)={n\choose m}\int_0^{\infty} (1-e^{-z})^m e^{-(n-m)z}\,{\Omega}_0({\rm d}z)
$$
and
$${q}_0(n,m)={n\choose m} {w}_0(n,m)+{w}_0(n,0)q(n,m).
$$
And when expressed  via binomial moments of $\omega$ this becomes
$$
{q}_0(n,m)= \frac{1}{\mu}{n\choose m}\left(\sum_{k=0}^m (-1)^{m-k}\frac{W(n-k)}{n-k}+\frac{w(n,m)}{n}
\right).
$$
\vskip0.5cm

\par {\large Remarks.} The relation between compositions ${C}_{0n}$ and $C_n$ is that they are
identically distributed given the size of the first part. 
The distribution of ${C}_n$ is of the form (\ref{pq}) with the first factor replaced by
${q}_0(n,n_1)$.

\par The distributional identity $(C_n)\stackrel{d}{=}(C_{0n})$ holds
iff $\Omega={\Omega}_0$, in which case $\Omega$ is an exponential distribution, ${\cal R}$ is a homogeneous 
Poisson  point process and therefore $C_n$ is governed by the ordered ESF.
This explains, to an extent, the role of ESF as a `central limit' because superposition of many 
rare renewal processes approaches the Poisson process.
\par For suitable choice of $\Omega$ the sums in the RHS of 
(\ref{ga}) or (\ref{ga+n}) become Cesaro or Euler averages. 
The LHS is easy to analyse but 
concluding directly from these relations about the
behaviour of $(a_n)$ is only possible
when $(a_n)$ is known to satify certain regularity conditions (the Tauberian conditions).
The direct approach   seems hard to realise 
because the regularity conditions are very sensitive
to the summability method. 

\vskip0.5cm

\par For $r_n\equiv 1$ the LHS of (\ref{ga+n}) is the expected number of parts of the stationary composition, 
which is equal to
$$1+\sum_{m=1}^{n-1} \frac{W(m)}{m\,\mu}$$ 
and, quite expectedly, is asymptotic to
$\mu^{-1} \log n$.

\vskip0.5cm
\par {\large Example.}
For ESF(1)   we have $W(n)=1-(n+1)^{-1}$ and $ \mu=1$ whence the expected number of parts is the
harmonic number 
$h_n=  1+(n-1)^{-1}+\ldots+3^{-1}+2^{-1}$ as is well-known \cite{ESF}.
\vskip0.5cm
\par We have seen that $g(n,m)\to g_0(n,m)$ for $n\to\infty$ and wish to obtain the asymptotics of (\ref{b}) by 
substituting $g_0$ instead of $g$. To this end we need  a stronger assumption on $\omega$
$$
\nu:=\int_0^1 (\log(1-x))^2\,\omega ({\rm d}x)<\infty\,
$$
which in terms  of  $\Omega$  means finiteness of the second moment, $\nu=\int_0^{\infty} x^2 \,\Omega({\rm d}x)$.

\vskip0.5cm

\begin{proposition} Suppose $\omega$ is not supported by a geometric sequence  and also $\nu<\infty$. 
Suppose $(r_n)$ is such that
$|r_n|<r_n'$ where $(r_n')$ is a decreasing sequence satisfying
$\sum r_n'/n<\infty$. Then for $(b_n)$ solving \rm{(\ref{b})} we have
$$
 \lim_{n\to\infty} b_n =  \frac{1}{\mu} \sum_{n=1}^\infty \frac{W(n) r_n}{n}$$
\end{proposition}
{\it Proof.} Given integer $J$ suppose $(r_n)$ is such that $r_n=0$ for $n<J$, is decreasing for $n\geq J$ and 
satisfies $\sum r_n/n<\infty.$ We wish to show that for the sequence 
$(b_n)$ solving (\ref{b}) with such $(r_n)$ there is a bound
\begin{equation}\label{bound}
\lim\,\sup\, b_n\,<\frac{1}{\mu}\sum_{n=1}^{\infty}\frac{ r_n}{n}+   \frac{r_J\nu}{\mu^2}\,.
\end{equation}
To this end we will make use of the renewal representation.

\par Recall that $Q_n$ collects reward $r_m$ if the chain visits state $m$. This occurs when ${\cal R}$
has at least one point between $E_{n-m}$ and $E_{n-m+1}$, in which case let us 
assign reward $r_m$ to  the rightmost such point
(equivalently, given $E_{n-m}$ is the leftmost point in a cluster, the 
 point of ${\cal R}$ in question is the left endpoint of the component interval$\subset [0,\infty]\setminus {\cal R}$ 
containing $E_{n-m}$).
let $U$ be the potential measure of
${\cal R}$, so that $U[0,z]$ is the expected cardinality of ${\cal R}\cap [0,z]$.
The total expected reward of $Q_n$ may be written as
   $$r_nW(n)+\int_0^{\infty} \Phi_n(z)U({\rm d}z)$$
where the first term stands for the reward at $0\in {\cal R}$,
which is only due in the event that  the first jump of the renewal process exceeds $E_1$,
and the integrand is
\begin{equation}\label{Phi}
\Phi_n(z)=\sum_{m=1}^{n}{n\choose m} e^{-zm}(1-e^{-z})^{n-m} r_m W(m).
\end{equation}
\par In the same manner, we associate the  rewards collected by  the stationary chain $Q_{0n}$ with the  separating 
 points of ${\cal R}_0$ and with $0$ (exceptional point, not in ${\cal R}_0$). The expected reward becomes
$$r_n \frac{W(n)}{n\mu}+     \int_0^{\infty} \Phi_n(z)U_0({\rm d}z)$$
where the first term stands for the event that $E_1<Z_0$ (i.e., $E_1$ falls to the left of ${\cal R}_0$).
This is, of course, yet another expression for the LHS of (\ref{ga+n}).

\par We modify now the reward processes for chains $Q_n$ and $Q_{0n}$ by deleting the first 
term (reward at $0$) and by replacing $\Phi_n$ with another function
$\widetilde{\Phi}_n$ defined via (\ref{Phi}) but with 
factors $W(m)$ deleted. Deleting the first term has no asymptotic effect because $r_n$ goes to $0$ as $n\to\infty$.
We also have  $\Phi_n(z)\leq \widetilde{\Phi}_n(z)$, thus
$\Phi_n(z)$ corresponds to a more generous reward structure,  with reward at $z$
 being $r_m$ if $E_{n-m}<z<E_{n-m+1}$ (thus there is no other contraint on $z$ except that $z\in {\cal R}$, respectively $z\in {\cal R}_0$).
The modified reward associated with ${\cal R}_0$ is the sum in the LHS of (\ref{bound}).
\par The function $\widetilde{\Phi}_n(z)$ is unimodal, with the unique maximum attained at $z^*$, which is the unique
positive solution of equation
$$
-r_J {n-1\choose J}+\sum_{m=J}^{n-1}(r_m-r_{m+1}){n-1\choose m}\left(\frac{e^{-z}}{1-e^{-z}}\right)^{m-J}=0
$$
(the uniqueness follows from the monotonicity of $(r_j), j>J$). For $n\to\infty$ Poisson approximation provides
asymptotics
$
n e^{-z^*}\to \zeta
$
where $\zeta$ is the unique positive root of transcendental equation
$$
-\frac{r_J}{J!}+\sum_{m=J}^{\infty}(r_m-r_{m+1})\frac{\zeta^{m-J}}{m!}=0.
$$
In the following argument it is only important that $z^*\to\infty$ as $n\to\infty$.

\par Because ${\cal R}_0$ is ${\cal R}$ shifted to the right, ${\cal R}_0=Z_0+{\cal R}$, there is a one-to-one correspondence between the sets
${\cal R}\,\cap \,](z^*-Z_0)_+,z^*]$ and ${\cal R}_0\,\cap \,]0,z^*]$. 
Furthermore, because
$\widehat{\Phi}_n$ is increasing
on $[0,z^*]$ the total (modified) reward of ${\cal R}_0$ over $]0,z^*]$ is larger than that of ${\cal R}$ on
$]0,(z^*-Z_0)_+]$. 
On the other hand, the expected reward of ${\cal R}$ on  
$](z^*-Z_0)_+\,,z^*]$ has an asymptotic bound $r_J\nu/(2\mu^2)$; indeed $r_J=\max r_j$ is an upper bound for the
instantaneous reward
and the potential $U\,](z^*-Z_0)_+\,,z^*]$ is asymptotic to 
$${\mathbb E}\frac{Z_0}{\mu}=
 \frac{1}{\mu}\int_0^{\infty} z\,\Omega[z,\infty]\,{\rm d}z=\frac{1}{2\mu}\int_0^{\infty}z^2\,\Omega({\rm d}z)
$$
as it follows from the two-term expansion in the renewal theorem, in the case  $\nu<\infty$
 (see
 \cite{Feller}, Section 4, Chapter XI). 
\par To the right of $z^*$ the relation is reversed, since the function $\widehat{\Phi}_n$ is decreasing. 
Shifting the origin to 
the  leftmost point of ${\cal R}_0\cap [z^*,\infty]$ enables to view ${\cal R}$ on the new scale as the range of a delayed renewal 
sequence. Thus the expected reward of ${\cal R}_0$ on $[z^*,\infty]$ is larger than that of ${\cal R}$, up to a term estimated by
$r_J\nu/(2\mu^2)$, exactly as above. Putting the two parts together shows  that the 
expected modified reward is bounded by 
the RHS of (\ref{bound}). The unmodified reward is smaller, hence (\ref{bound}).

\par Now suppose $(r_n)$ is decreasing and  satisfies $\sum r_n/n<\infty$.
 We split the sequence at $J$ and decompose it in two:
$r_n=r_n1_{\{n<J\}}+r_n1_{\{n\geq J\}}$.
Since recursion (\ref{b}) is linear, the decomposition forces   
the representation of solution as, say
 $b_n=b_n'+b_n''$. 
Applying the renewal theorem we get
$b_n'\to \mu^{-1}\sum_{n=1}^J r_nW(n)/n$.
As for the second part,
$\lim\,\sup\, b_n''$ is estimated with the help of (\ref{bound}) and approaches  zero when $J\to\infty$,
because both $r_J$ and the tail-sum of the series vanish.

\par For arbitrary sequence satisfying the condition of proposition splitting at 
$J$ yields one part converging to $\mu^{-1}\sum_{n=1}^J r_n W(n)/n$
and another part estimated by a solution with reward sequence decreasing for $n>J$, thus going to $0$ 
as $J$ grows.$\Box$

\vskip0.5cm

\par Now we are in a position to improve on the asymptotics of $a_n={\mathbb E}K_n$. It will not be
supererogatory to remind that asymptotic expansion of the harmonic number starts with $h_n=\log n+\gamma+O(n^{-1})$.

\begin{proposition}
Suppose $\omega$ is not supported by a geometric sequence, 
satisfies {\rm (\ref{llog})} and $\nu<\infty$. Then
$$a_n=\frac{\log n}{\mu}+\frac{\gamma}{\mu}+b+o(1)$$
where  $\gamma$ is the Euler constant and

$$b=\frac{1}{\mu}\int_0^1\log x\,\,\omega({\rm d}x)+\frac{\nu}{2\mu^2}.$$


\end{proposition}
{\it Proof.} Writing $a_n=\mu^{-1}h_n+b_n$, substituting this into (\ref{a}) 
and using the summation formula (\ref{firstsum})
we find that $(b_n)$ satisfies (\ref{br}) with 
$$r_n= 1-   \frac{1}{\mu W(n)}\int_0^1
\left(\frac{x^{}}{1}+
\frac{x^{2}}{2}+\ldots+\frac{x^n}{n}\right)\,\omega ({\rm d}x),$$
which can also be written  as 
$$
r_n W(n)=-\int_0^1 (1-x)^n \,\omega({\rm d}x)  +\frac{1}{\mu} \int_0^1 
\left(\frac{x^{n+1}}{n+1}+
\frac{x^{n+2}}{n+2}+\ldots\right)\,\omega ({\rm d}x).$$
Using monotone convergence and manipulating the series we find
$$\sum_{n=1}^{\infty}\frac{ r_n W(n)}{n\mu}=\frac{1}{\mu}\int_0^1\log x\,\omega\,({\rm d}x)+\frac{1}{2\mu^2}
\int_0^1(\log (1-x))^2\,\omega({\rm d}x)=b.$$
Since $W(n)\to 1$ and $r_nW(n)$ is the difference of two terms which decrease in $n$, application of Proposition 4 
 yields $b_n\to b$.
 $\Box$

\vskip0.5cm 

\par {\bf 6.} With no additional assumptions we will derive asymptotics of the variance
$v_n={\rm Var\,} K_n$. The key issue is the asymptotic evaluation of the inhomogeneous term of the recursion.

\vskip0.5cm

\begin{lemma} Under  assumptions of Proposition 5 the expectation $a_n={\mathbb E} K_n$ satisfies
$$
\lim_{n\to\infty}\left( 2a_n-1-a_n^2+\sum_{m=1}^{n-1} q(n,m) a_{n-m}^2\right) = \frac{\nu}{\mu^2}-1    \,.
$$
\end{lemma}

{\it Proof.} For $b_n,r_n$ bearing the same meaning  as in Proposition 5, we have
$$
b_n\to b,\quad W(n)\to 1,\quad b_n-\sum_{m=1}^{n}q(n,m) b_{n-m}=r_n,
$$
and integrating by parts yields
$$
r_nW(n)=-n \int_0^1 \omega\,[0,x]\,(1-x)^{n-1}{\rm d}x+\frac{1}{\mu}\int_0^1\frac{\omega\,[x,1]\, x^n}{1-x}\,{\rm d}x\,.
$$

\par Further useful estimates follow from the $n\to\infty$ asymptotics
\begin{equation}\label{log}
\int_0^1 x^n  \,\omega({\rm d}x)=o\left(\frac{1}{\log n}\right),\qquad
\int_0^1 (1-x)^n  \,\omega({\rm d}x)=o\left(\frac{1}{\log n}\right),\qquad r_n=o\left(\frac{1}{\log n}\right)\,.
\end{equation}
To justify the first  relation, observe that integrability and monotonicity of $\log\, (1-x)$ imply 
that $\omega\,[x,1]=o\left(|\log (1-x)|^{-1}\right)$ for $x\uparrow 1$ (in fact, the relation
is equivalent to the integrability).
Integrating by parts and using monotonicity we have
$$\int_0^1 x^n  \,\omega({\rm d}x)=
\int_0^1 n x^{n-1} \omega[x,1]\,{\rm d}x\,<{\rm const\,\cdot\,} 
\int_0^1 n x^{n-1} |\log x|^{-1} \,{\rm d}x
$$
and by a Tauberian argument this is
$o\left(|\log n|^{-1}\right)$.
The second  relation follows in the same way from
$$\int_0^1(1- x)^n  \,\omega({\rm d}x)=n\int_0^1(1- x)^{n-1} \omega[0,x]\,{\rm d}x\,
<{\rm const\,\cdot\,} 
\int_0^1 n (1-x)^{n-1} |\log (1- x)|^{-1} \,{\rm d}x
$$
and $\omega\,[0,x]=o\left(|\log \,x|^{-1}\right)$ for  $x\downarrow 0$.
And the third relation follows from the first two.

\par Substituting $a_n=\mu^{-1} h_n+b_n$ 
and grouping terms we have
$$2a_n-1-a_n^2+\sum_{m=1}^n q(n,m) a_{n-m}^2=
T_1+T_2+T_3-1$$
with three to-be-evaluated terms 
\begin{eqnarray*}
T_1& =&-b_n^2+\sum_{m=1}^{n}q(n,m) b_{n-m}^2\\
T_2& =& 2b_n-2b_n\frac{h_n}{\mu}-   \frac{2}{\mu}\sum_{m=1}^{n} q(n,m) b_{n-m}h_{n-m}\\
T_3& =& \frac{2h_n}{\mu}-\frac{h_n^2}{\mu^2}+\frac{1}{\mu^2}\sum_{m=1}^{n} q(n,m)h_{n-m}^2.
\end{eqnarray*}
\par From $b_n\to b$ it is obvious that $T_1\to 0$ as $n\to\infty$. To see that also  $T_2$ vanishes  write
$$b_{n-m}h_{n-m} =b_{n-m}h_{n}+(b_{n-m}-b)(h_{n-m}-h_n)+b(h_{n-m}-h_n)$$
then from (\ref{firstsum}) and (\ref{log})  we obtain
$$\sum_{m=1}^{n-1}q(n,m)(h_{n-m}-h_n)=
\frac{1}{W(n)}\int_0^1\left( h_nx^n-\sum_{j=1}^n 
\frac{x^j}{j}\right)\,\omega{\rm d}x\,,$$
hence by (\ref{log}) and Lemma 1
\begin{eqnarray*}
 b \,\sum_{m=1}^{n} q(n,m) (h_{n-m}-h_n)\to -b\mu \\
\sum_{m=1}^{n} q(n,m)(b_{n-m}-b)(h_{n-m}-h_n)\to 0 \\
\sum_{m=1}^{n} q(n,m)b_{n-m}h_{n}=(b_n-r_n)h_n=b_n h_n+o(1)\\
\end{eqnarray*}
which indeed implies $T_2\to 0$.

\par To evaluate $T_3$ we need a summation formula similar to (\ref{firstsum}), but this time  we should 
take a combinatorial analogue
of $\log^2$ in place of $\log$. To this end, introduce 
$$s_n=\sum_{1\leq i\leq j\leq n} \frac{1}{i\,j}$$
then there is a summation formula
\begin{equation}\label{secondsum}
\sum_{m=0}^{n-1} {n\choose m} x^m (1-x)^{n-m} s_{n-m}=s_n-\sum_{j=1}^n \frac{x^j}{j}(h_n-h_{j-1})
\end{equation}
where we recognise partial sum of the Taylor series 
$$
\frac{1}{2}(\log (1-x))^2=\sum_{j=1}^{\infty} \frac{x^j}{j}\,h_{j-1}\,.
$$
It follows that
$$
\sum_{m=1}^{n-1} q(n,m) s_{n-m}=s_n-\frac{1}{W(n)} \int_0^1\left( \sum_{j=1}^n \frac{x^j}{j}(h_n-h_{j-1})\right)\,
\omega({\rm d}x)
$$
and because $h_n^2$ differs from $2s_n$ by the partial sum of a converging series,
$$ h_n^2=2s_n-\sum_{j=1}^n \frac{1}{j^2}\,,$$
we conclude that 
\begin{eqnarray*}
\sum_{m=1}^{n-1} q(n,m) h^2_{n-m}& = & h_n^2-\frac{2}{W(n)}\int_0^1 \left( \sum_{j=1}^n \frac{x^j}{j}(h_n-h_{j-1})\right)+o(1)=\\
h_n^2-\frac{2\mu h_n}{W(n)} +\frac{\nu}{W(n)}+o(1)& = &
h_n^2-2\mu h_n + \nu+o(1)
\end{eqnarray*}
where we exploited monotone convergence and (\ref{log}). Now it is easily seen that $T_3\to \nu/\mu^2$.
\par Putting the terms together we arrive at $T_1+T_2+T_3-1\to 1-\nu/\mu^2$. $\Box$
\vskip0.5cm

\par {\large Remark.} The summation formula (\ref{secondsum}) implies an analogue of Lemma 1:
for arbitrary normalised weight $\omega$ 
the square of logarithm is
$L^2([0,1],\omega)$-approximable by its Bernstein polynomial.

\vskip0.5cm

\par Appealing to Corollary 1 we obtain the desired asymptotics of variance. 
Define $\sigma^2=\nu-\mu^2$, that is, $\sigma^2=\int_{[0,\infty]}(z-\mu)^2\,\Omega({\rm d}z)$
is the variance of distribution $\Omega$.

\vskip0.5cm

\begin{proposition} Under assumptions of Proposition 5
$${\rm Var\,} K_n\sim \frac{\sigma^2}{\mu^3}\,\log n\,.$$
\end{proposition}

\vskip0.5cm

\par {\bf 7.} We turn next to the central limit theorem for $K_n$. 
Neininger and R{\"u}schendorf \cite{nein} derived a general CLT 
for solutions of equations
like (\ref{K}). In our context,
the assumptions of their Theorem 2.1 are easily checked, with the only exception that their CLT requires
some expansion ${\rm  Var}\, K_n=\mu^{-1}\,\log n+O((\log n)^{1-\epsilon})$, which is not 
guaranteed by the integrability of $(\log (1-x))^2$ rather relies on integrability of a higher power
of  the logarithm.  
We shall see  that in our situation   no additional assumptions are necessary and the CLT follows 
by a simple comparison with the number of renewals.
\par  Given $n$, define a {\it cell} to be a component interval of
$[0,\infty]\setminus{\cal R}$ containing at least one $E_j$, $j\leq n$. Clearly, the total number of cells is $K_n$.
Let $L_n$ be the number of cells which have the left endpoint smaller $\log n$, and let $R_n$ be the number 
of renewals on $[0,\log n]$  (including $0$), that is 
$R_n=\# ({\cal R}\,\cap\, [0,\log n])$. It is an easy matter to see that $L_n\leq R_n$ and $L_n\leq K_n$.
Moreover, since the expected number of order statistics that exceed $\log n$ is 1, we have ${\mathbb E}\,(K_n-L_n)<1$.
\vskip0.5cm

\begin{proposition} Under assumptions of Proposition 5, 
\begin{equation}\label{norm}
\frac{K_n-\mu^{-1}\log n}{\sigma \mu^{-3/2}\,\log n}
\end{equation}
converges weakly to the standard normal random variable.
\end{proposition}
{\it Proof.} By \cite{Feller} (Section 5, Ch. XI), $R_n$ is asymptotically normal with expectation 
$\mu^{-1}\log n$ and variance $\sigma^2\mu^{-3}\log n$. Furthermore,
${\mathbb E}R_n=\mu^{-1}\log n+ \nu (2\mu^2)^{-1}+o(1)$ (\cite{Feller}, Equation (4.5)).
By asymptotics of moments (Propositions 5 and 6) and the above inequalities,
the $L^1$-distance between  any two of the three random
 variables $(K_n-a_n)v_n^{-1/2}$, $(L_n-a_n)v_n^{-1/2}$ and $(R_n-a_n)v_n^{-1/2}$ goes to zero.  
It follows that $L_n$ and $K_n$ are also asymptotically normal. $\Box$

\vskip0.5cm

\par In fact, the renewal theorem  taken together with a Poisson limit for the number of $E_j$'s
exceeding $\log n$ implies weak convergence of $K_n-L_n$. 
Asymptotics of the expectation involves the exponential  integral function
$$I(z)=\int_z^{\infty}e^{-y} y^{-1}\,{\rm d}y\,.$$

\begin{proposition} We have
$$\lim_{n\to\infty} {\mathbb E} (K_n-L_n)=\frac{\gamma}{\mu}+ \frac{1}{\mu}\int_0^1 \log x \,\omega({\rm d}x)
+\int_0^1 I(x)\,\omega ({\rm \, d}x)\,.
$$
\end{proposition}
{\it Proof.} 
Recalling (\ref{Phi}), using Poisson approximation and the renewal theorem, and changing the variable
 of integration for $\zeta=ne^{-z}$ we compute
\begin{eqnarray*}
{\mathbb E} (K_n-L_n)=
\frac{1}{\mu}\int_{\log n}^{\infty} \Phi_n (z)\,{\rm d}z+o(1)= \int_0^1 e^{-\zeta}\sum_{m=1}^{\infty} \frac{\zeta^m W(m)}{m!}
\frac{{\rm d}\zeta}{\mu\zeta}+o(1)=\\
 \int_{0}^1 \int_0^1  \frac{1-e^{-zx}}{\mu z} \,{\rm d}z\,\omega ({\rm d}x)+o(1)=\frac{\gamma}{\mu}+
\frac{1}{\mu}\int_0^1 \log x \,\omega({\rm d}x)+
\frac{\gamma}{\mu}\int_0^1 I(x)\,\omega({\rm d}x)+o(1)
\end{eqnarray*}
where we also used 
$$e^{-\zeta}\sum_{m=1}^{\infty}\frac{ \zeta^m (1-(1-x)^m)}{m!}=1-e^{-\zeta x}\,$$ 
and the well-known formula 
$$\int_0^x \frac{1-e^{-y}}{y}\,{\rm d} y=I(x)+\log x +\gamma\,.\Box$$

\par Now recalling
$${\mathbb E}K_n=\frac{\log n}{\mu}+\frac{\gamma}{\mu}+\frac{1}{\mu}\int_0^1 \log x\,\omega \,({\rm d}x)
+\frac{\nu}{2\mu^2}+o(1)$$
and comparing the expectations
$${\mathbb E}L_n=\frac{\log n}{\mu}+\frac{\nu}{2\mu^2}-\frac{1}{\mu}\int_0^1 I(x)\,\omega \,({\rm d}x)+o(1)\,,\qquad 
{\mathbb E} R_n=\frac{\log n}{\mu}+\frac{\nu}{2\mu^2}+o(1)$$
we not only confirm `by computation' the inequality ${\mathbb E}L_n\leq {\mathbb E}R_n$ ($=U[0,\log n]$) but also 
come to the conclusion that {\it the number of component intervals of $[0,\log n] \setminus {\cal R}$
which contain no $E_j$'s, $j\leq n$, 
remains bounded as $n\to\infty$}. This conclusion is in good accord with the general point taken in 
\cite{GP}, \cite{TSF} that the
composition $C_n$ is a proper combinatorial analogue of the regenerative set $\cal R$.$\Box$

\vspace{0.5cm}
\par {\bf Acknowledgement.} I would like to thank Andrey Levin for help with (\ref{secondsum}) 
and further Euler-summation formulas.

\vskip0.5cm 
\noindent
\\
gnedin@math.uu.nl

\end{document}